\documentclass{amsart}
\usepackage{amssymb}
\usepackage{color}

\newtheorem{theorem}{Theorem}
\newtheorem{lemma}{Lemma}
\newtheorem{proposition}{Proposition}
\newtheorem{remark}{Remark}
\newcommand{\pp}{\noindent {\bf Proof. }}
\begin{document}
\title{$\mathbb Z_2$-graded codimensions of unital algebras}

\author[D.D. Repov\v s and M.V. Zaicev]
{Du\v san D. Repov\v s and Mikhail V. Zaicev}

\address{Du\v san D. Repov\v s \\Faculty of Education, and
Faculty of  Mathematics and Physics, University of Ljubljana,
 Ljubljana, 1000, Slovenia}
\email{dusan.repovs@guest.arnes.si}

\address{Mikhail V. Zaicev \\Department of Algebra\\ Faculty of Mathematics and
Mechanics\\  Moscow State University \\ Moscow,119992, Russia}

\email{zaicevmv@mail.ru} 



\begin{abstract}
We study polynomial identities of nonassociative algebras constructed by using infinite binary
words and their combinatorial properties. Infinite periodic and Sturmian words were first applied
for constructing examples of algebras with arbitrary real  PI-exponent greater than one. Later we
used these algebras for confirmation of the conjecture that PI-exponent increases
precisely by one after adjoining an external unit to a given algebra. Here we prove the same result for
these algebras for graded identities and graded PI-exponent, provided that the grading group is cyclic
of order two.
\end{abstract}

\date{\today}

\maketitle\vskip 0.2in 

\section{Introduction}

We study numerical invariants of polynomial identities of algebras over a field of characteristic zero.
One of the most important characteristics of identities of an algebra $A$ is its codimension sequence
$\{c_n(A)\}$. In many cases this sequence is exponentially bounded and one can ask whether the limit
$$
exp(A)=\lim_{n\to\infty} \sqrt[n]{c_n(A)}
$$
exists. The answer is in general negative \cite{ERA}. Nevertheless, there is a wide class of algebras
where $exp(A)$ exists, for example, associative PI-algebras \cite{GZ1}, \cite{GZ2}, finite dimensional
Lie algebras \cite{GRZ1}, \cite{GRZ2}, \cite{Z}, finite dimensional Jordan and alternative algebras
\cite{GShZ}, and many others. Given an algebra $A$, one can consider an extension $A^{\#}$ of $A$, obtained
from $A$ by adjoining the external unit element. Then some natural questions arise: is the codimension sequence
$c_n(A^{\#})$ exponentially bounded, does $exp(A^{\#})$ exist, does there exist a relationship between $exp(A)$
and $exp(A^{\#})$?

It was first  mentioned in \cite{GZ3} that $exp(A^{\#})$ exists and is equal to either $exp(A)$ or $exp(A)+1$
for any associative PI-algebra $A$. One of the first examples of a non-associative algebra $A$ with
$exp(A^{\#})=exp(A)+1$ was found in \cite{Z-AL}. In the same paper it was conjectured  that for any
algebra $A$ either $exp(A^{\#})=exp(A)$ or $exp(A^{\#})=exp(A)+1$. In \cite{RZ-CA} this
conjecture was confirmed for a wide class of algebras associated with infinite  Sturmian words.
Similar results for certain kinds of Poisson algebras were found in \cite{Rats}. Note also that
$exp(A^{\#})=exp(A)$, whenever $A$ is a unital algebra \cite{BBZ}.

If $A$ is equipped with a group grading then one can also consider its graded identities and graded
codimensions $\{c_n^{gr}(A)\}$. In the present paper we begin to study connections between asymptotics of
$\{c_n(A)\}$ and $\{c_n^{gr}(A)\}$. We prove that for the class of algebras introduced in \cite{GMZAdv}
and associated with Sturmian words, graded PI-exponents exist, $exp^{gr}(A)=exp(A)$, and
$exp^{gr}(A^{\#})=exp^{gr}(A)+1$ for the most
natural $\mathbb Z_2$-grading. For all details concerning the polynomial identities and their numerical
invariants we refer to \cite{Dr}, \cite{GZbook}.

 \section{Preliminaries and main constructions}

Let $A$ be an algebra over a field $F$ of characteristic zero and let $F\{X\}$ be the absolutely free
algebra over $F$ with an infinite set of generators $X$. A polynomial $f=f(x_1,\ldots, x_n)$,
$x_1,\ldots,x_n\in X$, is called an identity of $A$ if $f(a_1,\ldots,a_m)=0$, whenever
$a_1,\ldots,a_n\in A$. The set $Id(A)$  of all identities of $A$ forms an ideal of $F\{X\}$. Denote
by $P_n$ the subspace of all multilinear polynomials on $x_1,\ldots, x_n$. Then $P_n\cap Id(A)$ is
the set of all multilinear identities of $A$ of degree $n$. It is well known that all identities of $A$
are completely defined by the family of subspaces $\{P_n\cap Id(A)\}, n=1,2,\ldots~$. An important numerical
invariant of identical relations of the algebra $A$ is the sequence of codimensions
$$
c_n(A)=\dim\frac{P_n}{P_n\cap Id(A)}.
$$
In the case of exponentially bounded growth of $\{c_n(A)\}$, one can define the lower and the upper PI-exponents
by setting
$$
\underline{exp}(A)=\liminf_{n\to\infty} \sqrt[n]{c_n(A)}, \quad
\overline{exp}(A)=\limsup_{n\to\infty} \sqrt[n]{c_n(A)}
$$
and the ordinary PI-exponent
$$
exp(A)=\lim_{n\to\infty} \sqrt[n]{c_n(A)},
$$
provided that $\underline{exp}(A)=\overline{exp}(A)$. A powerful tool for studying asymptotics
of codimensions is the representation theory of the symmetric group $S_n$. The group $S_n$ acts naturally on the space
$P_n$ of multilinear polynomials
$$
\sigma\, f(x_1,\ldots,x_n)=f(x_{\sigma(1)},\ldots,x_{\sigma(n)}).
$$
Under this action, the subspaces $P_n, P_n\cap Id(A)$ and the quotient
$$
P_n(A)=\frac{P_n}{P_n\cap Id(A)}
$$
become $S_n$-modules. Consider the $n$th cocharacter of $A$, that is the character of $P_n(A)$,
$\chi_n(A)=\chi(P_n(A))$, and its decomposition into irreducible components
\begin{equation}\label{eqp1}
\chi_n(A)=\sum_{\lambda\vdash n} m_\lambda\chi_\lambda,
\end{equation}
where $\chi_\lambda$ denotes the irreducible $S_n$-character, corresponding to the partition $\lambda$ of
$n$, and the integer $m_\lambda$ denotes its multiplicity in $\chi_n(A)$.

Denote by $d_\lambda=\deg \chi_\lambda$ the dimension of the irreducible $S_n$-module with the character
$\chi_\lambda$. It follows from (\ref{eqp1})  that
\begin{equation}\label{eqp2}
c_n(A)=\sum_{\lambda\vdash n} m_\lambda d_\lambda.
\end{equation}

Another important numerical characteristic of $Id(A)$ is its $n$th colength
$$
l_n(A)=\sum_{\lambda\vdash n} m_\lambda.
$$
In  many cases the sequence $l_n(A)$ is polynomially bounded while $d_\lambda$'s in (\ref{eqp2}) are
exponentially large.  This means that the asymptotics of $c_n(A)$ is actually defined by the maximal value
of $d_\lambda$ with $m_\lambda\ne 0$.

For group graded algebras, identical relations and corresponding numerical invariants can also be
considered. We restrict ourselves to the case of $\mathbb Z_2$-gradings. Consider the free algebra
$F\{X,Y\}$ with two independent sets of generators $X$ and $Y$. We can endow $F\{X,Y\}$ with a
$\mathbb Z_2$-grading, by setting $\deg x=0,\deg y=1$, for all $x\in X,y\in Y$, and extending this grading to
all monomials on $X\cup Y$. If $A=A_0\oplus A_1$ is a $\mathbb Z_2$-graded algebra over $F$ then a polynomial
$f(x_1,\ldots,x_k,y_1,\ldots,y_m)\in F\{X,Y,\}$ is a graded identity of $A$ if
$f(a_1,\ldots,a_k,b_1,\ldots,b_m)=0$, for all $a_1,\ldots,a_k\in A_0,b_1,\ldots,b_m\in A_1$.

The set of all graded identities of $A$ forms a homogeneous in $\mathbb Z_2$-grading ideal $Id^{gr}$ of
$F\{X,Y\}$. The intersection $P_{k,m}\cap Id^{gr}(A)$ consists of all multilinear graded identities
of degree $k$ on even variables and of degree $m$ on odd variables, where $P_{k,m}$ is the subspace
of all polynomials multilinear on $x_1,\ldots,x_k,y_1,\ldots,y_m$. As before, the symmetric groups $S_k, S_m$
act independently on even and odd variables and both $P_{k,m}$ and $P_{k,m}\cap Id^{gr}(A)$, and also
$$
P_{k,m}(A)=\frac{P_{k,m}}{P_{k,m}\cap Id^{gr}(A)}
$$
are $S_k\times S_m$-modules. One can decompose $P_{k,m}(A)$ into irreducible components and write
$$
\chi(P_{k,m}(A))=\sum_{{\lambda\vdash k\atop \mu\vdash m}} m_{\lambda,\mu}\chi_{\lambda,\mu},
$$
where $\chi_{\lambda,\mu}$ is the irreducible $S_k\times S_m$-character and $m_{\lambda,\mu}$ is its
multiplicity. It is well known that $\chi_{\lambda,\mu}=\chi_\lambda\otimes\chi\mu$ and that
$$
\deg \chi_{\lambda,\mu}=\deg\chi_\lambda\deg\chi_\mu=d_\lambda d_\mu.
$$

Partial codimensions and colengths are defined as follows:
$$
c_{k,m}(A)=\deg\chi(P_{k,m}(A))=\dim P_{k,m}(A),
$$
$$
l_{k,m}(A)=\sum_{{\lambda\vdash k\atop \mu\vdash m}} m_{\lambda,\mu}.
$$
Finally, the graded $n$th codimension and the colength of $A$ are equal to
$$
c_n^{gr}(A)=\sum_{k=0}^n {n\choose k} c_{k,n-k}(A)
$$
and
$$
l_n^{gr}(A)=\sum_{k=0}^n l_{k,n-k}(A),
$$
respectively.

Graded PI-exponents are defined similarly,
$$
\underline{exp}^{gr}(A)=\liminf_{n\to\infty} \sqrt[n]{c_n^{gr}(A)}, 
$$ 
$$
\overline{exp}^{gr}(A)=\limsup_{n\to\infty} \sqrt[n]{c_n^{gr}(A)}, 
$$ 
$$
exp(A)^{gr}=\lim_{n\to\infty} \sqrt[n]{c_n^{gr}(A)}.
$$

Generalizing (\ref{eqp2}), we get
\begin{equation}\label{eqp3}
c_{k,n-k}(A)=\sum_{{\lambda\vdash k\atop \mu\vdash n-k}} m_{\lambda,\mu}d_\lambda d_\mu.
\end{equation}

Graded and ordinary codimensions satisfy the relation
\begin{equation}\label{eqp4}
c_n(A)\le c_n^{gr}(A)
\end{equation}
(see \cite{G-R} or \cite{BDLaa}).

We will use the following auxiliary function for computing codimensions. Let $x_1,\ldots,x_d$ be
non-negative real numbers such that $x_1+\cdots+x_d=1$, $d\le 2$. Then
$$
\Phi(x_1,\ldots,x_d)=\frac{1}{x_1^{x_1}\cdots x_d^{x_d}}.
$$
If $d=2$ then we write
$$
\Phi(x_1,x_2)=\Phi(x)=\frac{1}{x^x(1-x)^{1-x}}
$$
instead of $\Phi(x_1,x_2)$, where $0\le x\le 1$.

\section{Sturmian words and Sturmian algebras}

In this section we recall the construction of algebras based on infinite binary words and
their combinatorial properties. First, let $K=k_1k_2\ldots$ be an infinite word with integers
$k_i\ge 2, i=1,2,\ldots$. Denote by $A(K)$ a non-associative algebra with the basis
$$
\{a,b, z_j^{(i)}\mid 1\le j\le k_i,\ i\ge 1\}
$$
and with the multiplication table given by
$$
z_1^{(i)}a=z_2^{(i)},\ldots,z_{k_i-1}^{(i)}a=z_{k_i}^{(i)},~
z_{k_i}^{(i)}b=z_1^{(i+1)},~i=1,2,\ldots.
$$
All other products are zero. Note that $A(K)$ is 2-step left nilpotent, that is
$x_1(x_2x_3)\equiv 0$ is an identity of $A(K)$. It allows us to omit brackets in all products
and write $x_1x_2x_3\cdots x_n$ instead of $(\cdots((x_1x_2)x_3)\cdots)x_n$, keeping in  mind that all
non-left normed products are zero. Algebras of this type are used intensively in the study of numerical
invariants of polynomial identities. For instance, in \cite{GMZAdv} the first examples of algebras with an
arbitrary exponential growth $\alpha^n, 1\le\alpha\in \mathbb R$, were presented. Examples of algebras
with an intermediate growth $n^{n^\beta}, 0<\beta<1$, were constructed in \cite{GMZ2}. Recently,
examples of commutative algebras with polynomial codimension growth $n^\alpha, 3<\alpha<4,$ were
presented in \cite{GMVZ}. Other important examples of abnormal codimension growth were constructed
in  \cite{GZForum}, \cite{ERA}.

In the present paper we study identities on algebras $A(K)$ of special kind. Let $m\ge 2$ be an integer and
let $w=w_1w_2\ldots$ be an infinite word in the alphabet $\{0;1\}$. We denote by $A(m,w)$ the algebra $A(K)$,
where $K$ is constructed as follows:
$$
k_i=m+w_i,~i=1,2,\ldots~~.
$$
Earlier, the algebras $A(m,w)$ have already been used for constructing a continuous family of unitary algebras
with non-integer PI-exponents and for confirmation of the conjecture that $exp(A^\#)=exp(A)+1$ (see \cite{RZ-CA}).

We recall some well known facts from the combinatorial theory of infinite words (see, for example, \cite{L}). Given
a binary word $w=w_1w_2\ldots$, the complexity $Comp_w$ of $w$ is the function
$Comp_w\colon\mathbb N\to \mathbb N$, where $Comp_w(n)$ is the number of distinct subwords of $w$ of length $n$.
It is easy to see that for a periodic word $w$ with period $T$, the complexity function is bounded,
$Comp_w(n)\le T$. Moreover, it is well known that $Comp_w(n)\ge n+1$ for any aperiodic $w$. An infinite
word $w$ is called {\it Sturmian} if $Comp_w(n)=n+1$ for all $n\ge 1$.

For a finite word  $x=x_1\ldots x_n$ in the alphabet $\{0;1\}$, the height $h(x)$ and the length $|x|$ are
defined as $h(x)=x_1+\cdots+x_n$ and $|x|=n$, respectively. Then the {\it slope} $\pi(x)$
is defined by
$$
\pi(x)=\frac{h(x)}{|x|}.
$$
One can extend this notion to certain infinite binary words. Namely, if the limit
$$
\pi(w)=\lim_{n\to\infty}\frac{h(w_1\ldots w_n)}{n}
$$
exists then $\pi(w)$ is called the slope of $w$. Clearly, the limit does not exist in general. Nevertheless,
for periodic and Sturmian words, the slope is well defined. In the next proposition we recall the basic
properties  which we will need in the sequel.
\begin{proposition}\label{P1}(\cite[Section 2.2]{L}).
Let $w$ be a Sturmian or periodic word. Then there exists a constant $C$ such that
\begin{itemize}
\item[1.] $|h(x)-h(y)|\le C$, for any finite subwords $x,y$ of $w$ with $|x|=|y|$;
\item[2.] the slope $\pi(w)$ of $w$ exists;
\item[3.] for any non-empty finite subword $u$ of $w$,
$$
|\pi(u)-\pi(w)|\le\frac{C}{|u|}; and
$$
\item[4.] for any real $\alpha\in (0;1)$, there exists a word $w$ with $\pi(w)=\alpha$ and $w$ is
Sturmian or periodic, according to whether $\alpha$ is irrational or rational, respectively.
\end{itemize}
\end{proposition}

We will use the following results.

\begin{theorem}\label{t1}(\cite[Theorem 5.1]{GMZAdv}
Let $w$ be a Sturmian or periodic word with the slope $0<\alpha<1$. If $m\le 2$ then for the algebra $A=A(m,w)$
the PI-exponent exists and $exp(A)=\Phi(\beta)$, where $\beta=\frac{1}{m+\alpha}$.
\end{theorem}
\hfill $\Box$

\begin{theorem}\label{t2}(\cite[Theorem 1]{RZ-MatSb})
Let $A=A(m,w)$, where $w$ is an infinite Sturmian or periodic word, and $m\ge 2$. Let  $A^\#$ be the algebra
obtained from $A$ by adjoining an external unit. Then PI-exponent of $A^\#$ exists and
$exp(A^\#)=exp(A)+1$.
\end{theorem}
\hfill $\Box$

\section{Gradings on Sturmian algebras}

The algebra $A=A(m,w)$ can be equipped by a $\mathbb Z_2$-grading in different ways. We begin our study with the most
natural case when generators of $A$ are homogeneous. The algebra $A$ is generated by the three elements
$z_1^{(1)}, a, b$. Each generator can be even or odd, so we have eight options. Clearly, if
$\deg z_1^{(1)}=\deg a=\deg b=0$ then the grading is trivial and all identities and codimensions are the same as
in the non-graded case. In the present paper we consider one of non-trivial cases when $z_1^{(1)}$ and $a$ are even,
whereas $b$ is odd. At the end we will discuss the difference between distinct gradings.

Throughout this section, let $A=A(m,w)$ be the algebra defined in the previous section, where $m\ge 2$ is an integer
and $w$ is an infinite periodic or Sturmian word. Then a $\mathbb Z_2$-grading $A=A_0\oplus A_1$ on $A$ is
uniquely defined by setting $\deg z_1^{(1)}=\deg a =0, \deg b=1$. First, we will give an upper bound for the graded
codimension $c_n^{gr}(A)$.

\begin{lemma}\label{l2}
Let $c_{k,n-k}(A)$ be the partial graded codimension of $A$. Then $c_{k,n-k}(A)\le 2n^2$ for all  large enough $n$.
\end{lemma}
\pp Consider a left-normed monomial $M=M(x_1,\ldots,x_k,y_1,\ldots,y_{n-k})$ on even
$x_1,\ldots,x_k$ and odd $y_1,\ldots,y_{n-k}$. Then $M=x_iu_1\cdots u_{n-1}$ or $M=y_iu_1\cdots u_{n-1}$, where
$u_1,\ldots, u_{n-1}$ are some $x_j$'s, $y_j$'s. Let for example,
$$
M=x_k\cdots x_{i_1}\cdots x_{i_{k-1}}\cdots~,~~\{i_1,\ldots,i_{k-1}\}=\{1,\ldots,k-1\}.
$$
Then $M\equiv M_0$ modulo the graded ideal $Id^{gr}(A)$, where $M_0=x_k\cdots x_1\cdots x_{k-1}\cdots$, since any
non-zero evaluation $\varphi$ of $M$ and $M_0$ in $A$ can be obtained only if $\varphi(x_k)=z_j^{(i)},
\varphi(x_1)=\cdots=\varphi(x_{k-1})=a$. Moreover, $\varphi(M)\ne 0$ if and only if the positions of
$x_{i_1},\ldots,x_{i_{k-1}}$ in $M$ are in 1-1 correspondence with the positions of symbol $0$ in the subword
$\bar w=w_{t+1}\cdots w_{t+n-1}$ of length $n-1$, where the integer $t$ can be computed from the condition
$z_j^{(i)}=z_1^{(1)}u_1\cdots u_t$ for proper $u_1,\ldots,u_t\in\{a,b\}$. Similarly, $y_1,\ldots,y_{n-k}$ in
$M$ can be ordered naturally. Since $Comp_w(n-1)=n$ for Sturmian word and $Comp_w$ is bounded in the
periodic case, we conclude that the number of subwords $\bar w$ corresponding to
monomials that do not vanish on $A$
does not exceed $kn\le n^2$ for  sufficiently large $n$. The same upper bound takes place for
monomials of the type $y_iu_1\cdots u_{n-1}$, and we have completed the proof.

\hfill $\Box$

\begin{lemma}\label{l3}
For any real number $\varepsilon>0$ there exists an integer $n_0$ such that conditions $n\ge n_0$,
$P_{k,n-k}(A)\ne 0$ imply the inequalities
\begin{equation}\label{eqg0}
\beta-\varepsilon\le\frac{n-k}{n}\le\beta+\varepsilon,
\end{equation}
where $\beta=\frac{1}{m+\alpha}$ and $\alpha=\pi(w)$ is the slope of the infinite word $w$ defining $A=A(m,w)$.
\end{lemma}

\pp
Any non-zero product of $n$ basis elements of $A$ has the form
\begin{equation}\label{eqg1}
z_j^{(i)}\underbrace{a\cdots a}_{s_0}b\underbrace{a\cdots a}_{s_1}b\cdots b
\underbrace{a\cdots a}_{s_r}b\underbrace{a\cdots a}_{s_{r+1}}= z_{1+s_{r+1}}^{(i+r+1)},
\end{equation}
where $0\le s_0,s_{r+1}\le m$,
\begin{equation}\label{eqg1a}
s_1=m+w_{i+1}-1,\ldots,s_r=m+w_{i+r}-1,
\end{equation}
$n=s_0+s_{r+1}+2+mr+w_{i+1}+\cdots+w_{i+r}$. The number of factors $b$ in this product is equal to
$r+1$. Moreover, (\ref{eqg1}) is the value of monomials from $P_{k,n-k}$ with $n-k=r+1$. Hence
\begin{equation}\label{eqg2}
\frac{n-k}{n}=\frac{r+1}{s_0+s_{r+1}+2+mr+w_{i+1}+\cdots+w_{i+r}}
\end{equation}
$$
=\frac{1+\frac{1}{r}}{m+\frac{s_0+s_{r+2}+2}{r}+\frac{w_{i+1}+\cdots+w_{i+r}}{r}}.
$$
Since $s_0+s_{r+1}+2\le 2m+2$ and $w_{i+1}+\cdots+w_{i+r}\le r$, it follows that
$$
r\ge \frac{n}{m+1}-2.
$$
In particular, $r\to\infty$ if $n\to\infty$. Moreover, the limit of $\frac{1}{r}(w_{i+1}+\cdots+w_{i+r})$,
as $r\to\infty$, is equal to $\alpha$, by Proposition \ref{P1}. It follows that the right hand side of
(\ref{eqg2}) goes to $\beta=\frac{1}{m+\alpha}$ as $n\to\infty$ and  hence (\ref{eqg0}) holds.

\hfill $\Box$

Lemmas \ref{l2} and \ref{l3} give an upper bound for graded codimensions of $A$.

\begin{lemma}\label{l4}
For any $0<\varepsilon\le\frac{1}{2}-\frac{1}{m+\alpha}$ there exists $n_0$ such that
$$
c_n^{gr}(A)\le 2n^3 \Phi(\beta+\varepsilon)^n
$$
for all $n\ge n_0$, where $\beta=\frac{1}{m+\alpha}$. In particular, $\overline{exp}^{gr}(A)\le\Phi(\beta)$.
\end{lemma}
\pp
By Lemmas \ref{l2} and \ref{l3} we have
$$
c_n^{gr}(A)\le \sum_{\beta-\varepsilon\le\frac{n-k}{n}\le \beta+\varepsilon}{n\choose k}
c_{k,n-k}(A) \le 2n^2 \sum_{\beta-\varepsilon\le\frac{n-k}{n}\le \beta+\varepsilon}
{n\choose k}.
$$
By Stirling's formula for factorials we have
$$
{n\choose k}\le n\frac{n^n}{k^k(n-k)^{n-k}}=n\Phi(\frac{n-k}{n})^n.
$$
Since $m\ge 2$ and $0<\alpha<1$, we have $\beta=\frac{1}{m+\alpha}<\frac{1}{2}$ and
$$
\max_{\beta-\varepsilon\le\frac{n-k}{n}} \Phi(\frac{n-k}{n})\le \Phi(\beta+\varepsilon)
$$
as soon as $\beta+\varepsilon<\frac{1}{2}$ and $n$ is sufficiently large. Hence
$$
c_n^{gr}(A)\le 2n^3\Phi(\beta+\varepsilon)^n,
$$
and we are done.
\hfill $\Box$

Now we are ready to prove main result of this section.

\begin{theorem}\label{t3}
Let $A=A(m,w)$ be the algebra defined by an integer $m\ge 2$ and by an infinite periodic or Sturmian
word $w$ with the slope $\pi(w)=\alpha$. Suppose that the decomposition $A=A_0\oplus A_1$ is a
$\mathbb Z_2$-grading of $A$ such that $a,z_1^{(1)}\in A_0, b\in A_1$. Then the graded PI-exponent
$exp^{gr}(A)$ exists and
$$
exp^{gr}(A)=exp(A)=\Phi(\frac{1}{m+\alpha}).
$$
\end{theorem}

\pp According to Lemma \ref{l4}, it is enough to show that $\underline{exp}^{gr}(A)
\ge\Phi(\beta)$, where $\beta=\frac{1}{m+\alpha}$.
Since $A$ is not nilpotent, there exists for any $n$, a non-zero product of the type
(\ref{eqg1}). In particular, given $n$, there exists $0\le k\le n$ such that $P_{k,n-k}\ne 0$. Then
$\frac{n-k}{n}\ge\beta-\varepsilon$ asymptotically for any fixed $\varepsilon>0$ by Lemma \ref{l3},
and by Stirling's formula we have
$$
c_n^{gr}(A)\ge{n\choose k}c_{k,n-k}(A) \ge{n\choose k}\ge\frac{1}{n^2}\frac{n^n}{k^k(n-k)^{n-k}}\ge
\frac{1}{n^2}\Phi(\beta-\varepsilon)^n.
$$It follows that $\underline{exp}^{gr}(A)\ge\Phi(\beta)$, and thus the proof has been completed.
\hfill$\Box$

\section{Algebras with adjoint unit}

In this section we study codimensions of algebras with an external unit. Given an algebra $B$, we denote by
$B^\#$ the algebra obtained by adjoining the external unit to $B$. Note that if $C=\oplus_{g\in G}C_g$ is a
$G$-graded algebra with unit $1$ then $1$ is a homogeneous element and $1\in C_e$, where $e\in G$ is the
identity element of the group  $G$. Therefore in the case of a $\mathbb Z_2$-graded algebra $B$, its extension
$B^\#=B\oplus 1$ has a unique $\mathbb Z_2$-grading $B^\#=B^\#_0\oplus B^\#_1$, where $B$ is a homogeneous
subalgebra of $B^\#$, namely, $B^\#_0=B_0\oplus 1, B^\#_1=B_1$.

First, let $A_0\oplus A_1$ be an arbitrary $\mathbb Z_2$-graded algebra. Denote by $R\{X,Y\}$ the relatively free
$\mathbb Z_2$-graded algebra of the variety $var^{gr}(A)$ of graded algebras generated by $A$ with two
infinite sets $X$ and $Y$ of even and odd generators, respectively. That is, $R\{X,Y\}=F\{X,Y\}/Id^{gr}(A)$.
Consider a partial $(k,n-k)$-cocharacter of $A$,
\begin{equation}\label{equ1}
\chi_{k,n-k}(A)=\chi(P_{k,n-k}(A))=\sum_{{\lambda\vdash k\atop \mu\vdash n-k}} m_{\lambda,\mu}\chi_{\lambda,\mu}.
\end{equation}

In order to bound the multiplicities $m_{\lambda,\mu}$ in (\ref{equ1}) we denote by $R_{d_0,d_1}^{k,n-k}(A)$
the subspace of polynomials on $X_{d_0}=\{x_1,\ldots,x_{d_0}\}$, $Y_{d_1}=\{y_1,\ldots,y_{d_1}\}$ in
$R\{X,Y\}$ of total degree $k$ on $X_{d_0}$ and total degree $n-k$ on $Y_{d_1}$. The same argument as in
\cite{VMGU}  gives us the next lemma.

Recall that the height $h(\lambda)$ of a partition $\lambda=(\lambda_1,\ldots,\lambda_t)$ is the number $t$
of its parts.

\begin{lemma}\label{l4a}
Let $m_{\lambda,\mu},\lambda\vdash k,\mu\vdash n-k$, be the multiplicity from (\ref{equ1}) with $h(\lambda)\le d_0,
h(\mu)\le d_1$. Then
$$
m_{\lambda,\mu}\le\dim R_{d_0,d_1}^{k,n-k}(A).
$$
\end{lemma}
\hfill $\Box$

Now, let $A=A(m,w)$ be the algebra defined by an infinite binary word $w=w_1w_2\ldots~$. The following lemma holds for
any $w$, not necessarily Sturmian or periodic.

\begin{lemma}\label{l5}
Suppose that $A$ is $\mathbb Z_2$-graded and that $\deg z_1^{(1)}=\deg a=0, \deg b =1$. Then
$$
\dim R_{d_0,d_1}^{k,n-k}(A)\le d_0m^2k Comp_w(n-k).
$$
\end{lemma}
\pp Denote by $W$ the span of all monomials
\begin{equation}\label{equ2}
x_0u_1\cdots u_{n-1}
\end{equation}
in $R\{X,Y\}$, where $u_1,\ldots,u_{n-1}\in X_{d_0}\cup Y_{d_1}$, $u_{i_1},\ldots,u_{i_{k-1}}\in X_{d_0}$ for some
$$i_1,\ldots,i_{k-1}\in\{1,\ldots,n-1\},$$ 
while $u_j\in Y_{d_1}$, provided that $j\ne i_1,\ldots,i_{k-1}$. Clearly,
$\dim R_{d_0,d_1}^{k,n-k}(A)\le d_0\dim W$.

Let $f=f(x_0,\ldots,x_{d_0},y_1,\ldots,y_{d_1})\in F\{X,Y\}$ be a linear combination of monomials of
the same type as (\ref{equ2}). Then $f\equiv 0$ is an identity  of $A$ if and only if $\sigma(f)=0$ for any
homomorphism $\sigma\colon F\{X,Y\}\to A$ such that
$$
\sigma(x_0)= z_j^{(i)}, \sigma (x_s)=a,\sigma(y_s)=b.
$$
Hence $\dim W$ does not exceed the codimension of the intersection of all ${\rm Ker}~\sigma$ in
$F_{d_0,d_1}^{k,n-k}$,   where
$F_{d_0,d_1}^{k,n-k}$ is a subspace of  $F\{X,Y\}$ defined similarly as $R_{d_0,d_1}^{k,n-k}(A)$. 
Consider the family
of graded homomorphisms $\varphi_{ij}\colon F\{X,Y\}\to A$ such that 
$$\varphi_{ij}(x_0)=z_j^{(i)}, \varphi_{ij}(x_s)=a,
\varphi_{ij}(y_s)=b,$$ for all $x_s\in X,y_s\in Y$. Then either $\varphi_{ij}(x_0u_1\cdots u_{n-1})=0$ or
$\varphi_{ij}(x_0u_1\cdots u_{n-1})=z_{1+s_{r+1}}^{(i+r+1)}$, the element from (\ref{eqg1}). The latter equality
takes place if and only if $s_0=m-1-j+w_i, 0\le s_{r+1}\le m-1+w_{i+r+1}, n=s_0+s_{r+1}+2+mr+w_{i+1}+\cdots+w_{i+r}$,
relations (\ref{eqg1a}) hold and all $x_1,\ldots,x_{d_0}$ stay on ``correct'' positions among $u_1,\ldots,u_{n-1}$,
according to the word $w$. In particular, ${\rm codim}~{\rm Ker}~\varphi_{ij}$ in $F_{d_0,d_1}^{k,n-k}$ is less than or equal
to one. Moreover, ${\rm Ker}~\varphi_{ij}={\rm Ker}~\varphi_{i'j}$ if the subwords $w_{i+1}\cdots w_{i+r+1}$, $w_{i'+1}\cdots w_{i'+r+1}$
coincide. It follows that the codimension of $\bigcap {\rm Ker}~\varphi_{ij}$ in $F_{d_0,d_1}^{k,n-k}$ is at most
$m^2Comp_w(r+1)=m^2 Comp_w(n-k)$. Since $\dim W$ is equal to ${\rm codim}~\bigcap{\rm Ker}~\varphi_{ij}$,
we have completed the proof of the lemma.
\hfill $\Box$

Next, we will find an upper bound for $\dim R_{d_0,d_1}^{k,n-k}(B^\#)$ in terms of
$\dim R_{d_0,d_1}^{k,n-k}(B)$ if $B$ is a $\mathbb Z_2$-graded algebra.

\begin{lemma}\label{l6}
Given a $\mathbb Z_2$-graded algebra $B$, suppose that $\dim R_{d_0,d_1}^{k,n-k}(B)
\le \theta k(n-k)^T$ for all $0\le k \le n$ and for some constant $\theta$. Then
$$
\dim R_{d_0,d_1}^{k,n-k}(B^\#)\le \theta(k+1)^{d_0+2}(n-k+1)^{T+d_1}.
$$
\end{lemma}
\pp Note that a multihomogeneous polynomial $f(x_1,\ldots,x_{d_0},y_1,\ldots,y_{d_1})$
is a graded identity of $B^\#$ if and only if all multihomogeneous on
$x_1,\ldots,x_{d_0},y_1,\ldots,y_{d_1}$ components of
$f(1+x_1,\ldots,1+x_{d_0},y_1,\ldots,y_{d_1})\in F\{X,Y\}^\#$ are identities of $B$. The total
number of such components does not exceed $(k+1)^{d_0}(n-k+1)^{d_1}$, provided that the
degree on $\{x_1,\ldots,x_{d_0}\}$ is at most $k$ and the degree on $\{y_1,\ldots,y_{d_1}\}$
is equal to $n-k$.

Let $f_1,\ldots,f_N\in F_{d_0,d_1}^{k,n-k}$. Consider the linear combination $f=\lambda_1 f_1+
\cdots+ \lambda_Nf_N$ with unknown coefficients $\lambda_1,\ldots,\lambda_N$. Any multihomogeneous
component $g=g(x_1,\ldots,x_{d_0},y_1,\ldots,y_{d_1})$ of $f(1+x_1,\ldots,1+x_{d_0},y_1,\ldots,y_{d_1})$
gives us at most
$$
\dim R_{d_0,d_1}^{j,n-k}(B)\le\theta j(n-k)^T
$$
linear equations on $\lambda_1,\ldots,\lambda_N$, provided that $g\equiv 0$ is an identity of $B$
and the degree of $g$ on $x_1,\ldots,x_{d_0}$ is equal to $j$. Hence $f\equiv 0$ is an identity
of $B^\#$ if $\lambda_1,\ldots,\lambda_N$ satisfy no more than $\widetilde N$ linear equations,
where
$$
\widetilde N=(k+1)^{d_0}(n-k)^{d_1}\theta(n-k)^T \sum_{j=0}^k j.
$$
Note that
\begin{equation}\label{equ3}
 \widetilde N \le\theta (k+1)^{d_0+2}(n-k+1)^{T+d_1}.
\end{equation}
Therefore if $N$ is greater than the right hand side of (\ref{equ3}) then $f_1,\ldots,f_N$
are linearly dependent modulo $Id^{gr}(B^{\#})$ and we have completed the proof.
\hfill $\Box$

Now we are ready to get an upper bound for graded colength of $A^\#$.

\begin{lemma}\label{l7}
Let $A=A(m,w)$, where $m\ge 2$ is an integer and $w$ is an infinite periodic or Sturmian word. Then
$$
l_{k,n-k}(A^\#)\le 3m^2(k+1)^8(n-k+1)^6
$$
and
$$
l_n^{gr}(A^\#)\le 3m^2(n+1)^{15}.
$$
\end{lemma}
\pp
Consider a partial cocharacter of $A^\#$
\begin{equation}\label{equ4}
\chi_{k,n-k}(A^\#)=\sum_{{\lambda\vdash k\atop \mu\vdash n-k }}
m_{\lambda,\mu} \chi_{\lambda,\mu}.
\end{equation}
The linear subspace $I=Span\{z_j^{(i)}\vert~ i,j\ge 1 \}$ forms a homogeneous ideal of $A^\#$
with zero multiplication and
$$
\dim\left(A^\#/I \right)_0=2,\quad \dim\left(A^\#/I \right)_1=1.
$$
Hence any multilinear polynomial alternating on $4$ even variables or on $3$ odd variables is
an identity of $A^\#$. Standard argument implies that $m_{\lambda,\mu}\ne 0$ in (\ref{equ4})
only if $h(\lambda)\le 3, h(\mu)\le 2$. By Lemma \ref{l5} we have
$$
\dim R_{3,2}^{k,n-k}(A)\le 3m^2k(n-k+1)\le 3m^2k(n-k)^2.
$$
Then by Lemmas \ref{l4a} and \ref{l6},
$$
m_{\lambda,\mu}\le\dim R_{3,2}^{k,n-k}(A^\#)\le 3m^2k(k+1)^5(n-k+1)^4.
$$
The number of summands on the right hand side of (\ref{equ4}) is not greater than $k^3(n-k)^2$,
hence
$$
l_{k,n-k}(A^\#)\le 3m^2k(k+1)^8(n-k+1)^6
$$
and
$l_n^{gr}(A^\#)\le 3m^2(n+1)^{15}$.
\hfill $\Box$

Now we specify necessary conditions for inequality $m_{\lambda,\mu}\ne 0$ in (\ref{equ4}).

\begin{lemma}\label{l8}
Let $A=A(m,w)$ and suppose that $w$ is a Sturmian or periodic word with the slope $\alpha$. Suppose that
$m_{\lambda,\mu}\ne 0$ in (\ref{equ4}), where $\lambda\vdash k,\mu\vdash n-k$. Then
$\lambda$ and $\mu$ satisfy the following conditions:
\begin{itemize}
\item[1.] $\lambda=(\lambda_1,\lambda_2,\lambda_3)$ with $\lambda_3\le 1$;
\item[2.] $\mu=(\mu_1,\mu_2)$ with $\mu_2\le 1$;
\item[3.] $\lambda_1+\lambda_2+\mu_1=n$ or $n-1$; and
\item[4.] for any $0<\varepsilon<\frac{1}{2}-\beta$ there exists an integer $n_0$ such that
$$
\mu_1\le\frac{\beta+\varepsilon}{1-\beta-\varepsilon}\lambda_1
$$
for all $n\ge n_0$, where $\beta=\frac{1}{m+\alpha}$.
\end{itemize}
\end{lemma}
\pp
Any multilinear polynomial $f$ containing an alternating set of order $4$ on even variables
vanishes on $A^\#$. If $A^\#$ contains two alternating sets of order $3$ on even variables
then also $f\in Id^{gr}(A^\#)$. From the structure of essential idempotents of group ring $FS_k$
it follows that $\lambda_4=0,\lambda_3\le 1$. This proves 1. Similar argument gives us 2.

Let $\lambda=(\lambda_1,\lambda_2,1),\mu=(\mu_1,1)$. If $M$ is an irreducible
$S_k\times S_{n-k}$-submodule of $P_{k,n-k}\subset F\{X,Y\}$ with the character
$\chi(M)=\chi_{\lambda,\mu}$ then $M$ is generated by a polynomial $f(x_1,\ldots,x_k,
y_1,\ldots,y_{n-k})$ alternating on $x_1,x_2,x_3$ and on $y_1,y_2$. If we evaluate
$x_1,x_2,x_3$ on $\{1,a\}$ and $y_1,y_2$ on $\{b\}$ then we get zero. Otherwise,  $x_1,x_2,x_3$
should be equal to $1,a, z_j^{(1)}$, and $y_1,y_2$ should be equal to $b, z_s^{(r)}$. In this
case the value is also zero. Hence $\lambda_3+\mu_2\le 1$, and we obtain 3.

Let us prove 4. If $m_{\lambda,\mu}\ne 0$ then there exists a polynomial 
$$f=f(x_1,\ldots, x_k,
y_1,\ldots y_{n-k})$$ which generates an irreducible $S_k\times S_{n-k}$-module with the
character $\chi_{\lambda,\mu}$ and an evaluation $\varphi\colon X\cup Y\to A^\#$ such that
$\varphi(f)\ne 0$. Moreover, $f$ contains $\lambda_2$ disjoint alternating sets of $x's$ of order $2$.
The set of values $\{\varphi(x_1),\ldots,\varphi(x_k)\}$ contains $\lambda_1\ge p\ge \lambda_2-1$
elements $a$, at most one element $z_j^{(i)}$ and $k-p$ or $k-p-1$ units. The set
$\{\varphi(y_1),\ldots,\varphi(y_{n-k})\}$ contains at most one odd $z_j^{(i)}$ and $q=\mu_1-1$ or
$\mu_1$ elements $b$.

Furthermore, there exists a non-zero product
$$
g=z_j^{(i)}a\cdots aba\cdots ab\cdots ba\cdots a
$$
which is equal (up to a scalar factor) to $\varphi(f)$. Denote by $p=\deg_ag, q=\deg_bg$ the numbers
of entries of $a$ and $b$ in $g$, respectively. If the total degree $N=\deg g=1+p+q$ increases then
there exists a corellation between the growth of $p$ and $q$ (provided that $g\ne 0$). Namely,
$$
\lim_{n\to\infty}\frac{q}{q+p}=\beta=\frac{1}{m+\alpha},
$$
hence
$$
\lim_{n\to\infty}\frac{q}{p}=\frac{\beta}{1-\beta}.
$$
It follows that there exists $r$ such that for any $q\ge r+1$ (and for corresponding $p$) we have
$$
\frac{q}{p}\le\frac{\beta+\varepsilon/2}{1-(\beta+\varepsilon/2)}\qquad{\rm and}\qquad
\frac{1}{p}+\frac{\beta+\varepsilon/2}{1-(\beta+\varepsilon/2)}\le
\frac{\beta+\varepsilon}{1-(\beta+\varepsilon)}.
$$
Since $\mu_1-1\le q$ and  $p\le\lambda_1$, we get
$$
\frac{\mu_1-1}{\mu_1}\le\frac{q}{p}\le \frac{\beta+\varepsilon/2}{1-(\beta+\varepsilon/2)}
$$
and
$$
\frac{\mu_1}{\lambda_1}\le\frac{q}{p}+\frac{1}{\lambda_1}\le\frac{q}{p}+\frac{1}{p}\le
\frac{\beta+\varepsilon}{1-(\beta+\varepsilon)},
$$
provided that $\mu_1\ge r$.

On the other hand, if $\mu_1< r$ then
$$
\frac{n}{\mu_1}>\frac{n}{r}\qquad{\rm and}\qquad \frac{\lambda_1}{\mu_1}>\frac{n}{2r}-1
$$
since $2\lambda_1+\mu_1\ge\lambda_1+\lambda_2+\mu_1\ge n-1$. Denote for short
$\gamma=\frac{\beta+\varepsilon}{1-(\beta+\varepsilon)}$. Then for all $n\ge
\frac{2(\gamma+1)}{\gamma}r$ we have
$$
\frac{n}{2r}\ge\frac{1}{\gamma}+1, \quad\frac{\lambda_1}{\mu_1}>\frac{1}{\gamma}.
$$
This proves 4.
\hfill $\Box$

In order to get an upper bound for graded codimensions we need some properties of the function
$\Phi(x_1,\ldots,x_d)$ introduced in Section 2. Recall that $\Phi(x_1,x_2,x_3)=x_1^{-x_1}x_2^{-x_2}x_3^{-x_3}$,
where $0\le x_1,x_2,x_3\in \mathbb R$, $x_1+x_2+x_3=1$.

\begin{lemma}\label{l9}
Let $x_3=\gamma x_2$ for a fixed coefficient $\gamma$. Then
$$
\max\Phi(x_1,x_2,x_3)=\frac{1+\gamma}{\gamma^{\gamma/(\gamma+1)}}+1.
$$
\end{lemma}

\pp Denote $x=x_1$. Then the relations $x_1+x_2+x_3=1$, $x_3=\gamma x_2$ imply
$$
x_2=\frac{1-x}{1+\gamma},\quad x_3=\frac{\gamma}{1+\gamma}(1-x).
$$
Denote also $\Phi(x_1,x_2,x_3)=f(x)$. Then
$$
f^{-1}(x)=x^x \left(\frac{1-x}{1+\gamma}\right)^\frac{1-x}{1+\gamma}
 \left(\frac{\gamma(1-x)}{1+\gamma}\right)^\frac{\gamma(1-x)}{1+\gamma}
$$
and
$$
g(x)=\ln f^{-1}(x)=x\ln x +\frac{1-x}{1+\gamma}\ln\frac{1-x}{1+\gamma}+
\frac{\gamma(1-x)}{1+\gamma}\ln\frac{\gamma(1-x)}{1+\gamma}.
$$
Hence
$$
g'(x)=\ln\frac{x}{(\frac{1-x}{1+\gamma})^\frac{1}{1+\gamma}
(\frac{\gamma(1-x)}{1+\gamma})^\frac{\gamma}{1+\gamma} }
$$
and $g'({\widetilde x})=0$ only if
$$
\widetilde x =
\left(\frac{1-\widetilde x}{1+\gamma}\right)^\frac{1}{1+\gamma}
\left(\frac{\gamma(1-\widetilde x)}{1+\gamma}\right)^\frac{\gamma}{1+\gamma}
\gamma^\frac{\gamma}{1+\gamma} = (1-\widetilde x)\rho,
$$
where
$$
\rho=\frac{\gamma^\frac{\gamma}{1+\gamma}}{1+\gamma}.
$$
That is,
$$
\widetilde x=\frac{\rho}{1+\rho}.
$$

Since
$$
g'(x)=\ln\frac{x}{1-x}+const
$$
on the interval $(0;1)$, we see that $f^{-1}(x)$ has a local minimum in $\widetilde x$. Direct computations
show that $f^{-1}(\widetilde x)=\widetilde x$ and
$$
\max\Phi =f(\widetilde x)=\frac{1}{\widetilde x}=1+\frac{1}{\theta}
=\frac{1+\gamma}{\gamma^\frac{\gamma}{1+\gamma}}+1.
$$
\hfill $\Box$

Now we are ready to compute the required upper bound for the upper PI-exponent.

\begin{remark}\label{r1}
If we denote $\frac{\gamma}{\gamma+1}$ by $\theta$ then $\gamma=\frac{\theta}{1-\theta}$.
In this case direct computations give us
$$
\frac{1+\gamma}{\gamma^\frac{\gamma}{1+\gamma}} = \frac{1}{\theta^\theta (1-\theta)^{1-\theta}}
=\Phi(\theta).
$$
Moreover, if $\gamma_1<\gamma_2\le 1$ then $\theta_1<\theta_2$ and $\Phi(\theta_1)<\Phi(\theta_2)$.
\end{remark}

\begin{lemma}\label{l10}
$$
\overline{exp}^{gr}(A^\#)\le exp^{gr}(A)+1.
$$
\end{lemma}

\pp
By (\ref{eqp3})
\begin{equation}\label{equ5}
c_n^{gr}(A^\#)\le l_n^{gr}(A^\#)\sum_k\sum_{{\lambda\vdash k,\mu\vdash n-k\atop m_{\lambda,\mu\ne 0}}}
{n\choose k} d_\lambda d_\mu.
\end{equation}

First estimate a fixed summand ${n\choose k}d_\lambda d_\mu$ provided that $m_{\lambda,\mu}\ne0$.
By Lemma \ref{l8}, we have $\lambda=(\lambda_1,\lambda_2,\lambda_3)$, $\lambda_3\le 1$, $\mu=(\mu_1,\mu_2)$,
$\mu_2\le 1$. By the Hook formula for degree of an irreducible representation,
$$
d_\lambda\le \frac{k!}{\lambda_1!\lambda_2!},\quad d_\mu\le n-k-1\le n.
$$
Since $n-k=\mu_1$ or $\mu_1-1$, we have $n-k+1\ge \mu_1$ and $n(n-k)!\ge \mu_1!$. Also, $n-2\le\lambda_1+\lambda_2+\mu_1\le n$,
that is $n!\le(\lambda_1+\lambda_2+\mu_1)!(n+2)^2$. Therefore
\begin{equation}\label{equ6}
{n\choose k}d_\lambda d_\mu\le \frac{n!}{k!(n-k)!}\cdot\frac{k!}{\lambda_1!\lambda_2!}\cdot n
\le (n+2)^4\frac{(\lambda_1+\lambda_2+\mu_1)!}{\lambda_1!\lambda_2!\mu_1!}.
\end{equation}
By the Stirling's formula
\begin{equation}\label{equ7}
\frac{(\lambda_1+\lambda_2+\mu_1)!}{\lambda_1!\lambda_2!\mu_1!}\le
n\frac{(\lambda_1+\lambda_2+\mu_1)^{\lambda_1+\lambda_2+\mu_1}}{\lambda_1^{\lambda_1}\lambda_2^{\lambda_2}\mu_1^{\mu_1}}
\le n\Phi(x_1,x_2,x_3)^n,
\end{equation}
where
$$
x_1=\frac{\lambda_2}{\lambda_1+\lambda_2+\mu_1},\quad
x_2=\frac{\lambda_1}{\lambda_1+\lambda_2+\mu_1},\quad
x_3=\frac{\mu_1}{\lambda_1+\lambda_2+\mu_1}.
$$
Denote $\mu_1/\lambda_1=\gamma$. Then $x_3=\gamma x_2$, and by Lemma \ref{l9} and
Remark \ref{r1},
$$
\Phi(x_1,x_2,x_3)\le\Phi(\theta)+1,
$$
where $\theta=\frac{\gamma}{\gamma+1}$. Fix an arbitrary small $\varepsilon>0$. We can assume that
$$
\frac{\beta+\varepsilon}{1-\beta-\varepsilon}< 1.
$$
Then by Lemma \ref{l8} we get that $\gamma<1$, $\theta$ is an increasing function of $\gamma$ on interval
$(0;1)$ and $\theta\le\frac{1}{2}$. Hence $\Phi(\theta)$ is also an increasing function of $\gamma$ and
$$
\Phi(x_1,x_2,x_3)\le 1+\Phi(\beta+\varepsilon)
$$
for all sufficiently large $n$. Applying (\ref{equ5}), (\ref{equ6}), (\ref{equ7}) and Lemma \ref{l7}
and taking into account that the number of partitions $\lambda\vdash k$ with $h(\lambda)\le 3,
\lambda_3\le 1$ is not greater than $n$, we obtain
$$
c_n^{gr}(A^\#)\le 6m^2(n+2)^{22}(1+\Phi(\beta+\varepsilon))^n,
$$
from which it follows that
$$
\overline{exp}^{gr}(A^\#)\le 1+\Phi(\beta)= 1+exp^{gr}(A).
$$
\hfill $\Box$

\begin{theorem}\label{t4}
Let $A=A(m,w)=A_0\oplus A_1$ be the algebra defined by an integer $m\ge 2$ and by Sturmian or periodic
word $w$ equipped with a ${\mathbb Z}_2$-grading, where generators $z_1^{(1)}$ and $a$ are even whereas $b$ is odd.
Let  $A^\#$ be obtained from $A$ by adjoining the external unit. Then its graded PI-exponent exists and
$$
exp^{gr}(A^\#)=1+exp^{gr}(A).
$$
\end{theorem}
\pp By \cite[Theorem 1]{RZ-CA}, $exp(A^\#)=exp(A)+1$. Hence $\underline{exp}^{gr}(A^\#)\ge
exp(A^\#)=exp(A)+1=exp^{gr}(A)+1$ by (\ref{eqp4}) and Theorem \ref{t3}. Now our statement follows from
Lemma \ref{l10}.
\hfill $\Box$

In conclusion, we discuss other ${\mathbb Z}_2$-gradings on $A=A(m,w)$. In the proof of Theorems
\ref{t3} and \ref{t4} we have never used the fact that $\deg z_1^{(1)}=0$. Hence the same results hold for
graded codimensions if $\deg z_1^{(1)}=\deg b=1,\deg a=0$. 

By slightly modifying arguments, one can prove
Theorems \ref{t3} and \ref{t4}, provided that $\deg a=1,\deg b=0$. Finally, if $\deg a=\deg b=1$ then the argument
is similar to that of \cite{GMZAdv}, \cite{RZ-CA}. Therefore we can generalize Theorems \ref{t3} and
\ref{t4} as follows.

\begin{theorem}\label{t5}
Let $A=A(m,w)$ be the algebra defined by an integer $m\ge 2$ and by an infinite periodic or Sturmian
word $w$ with the slope $\pi(w)=\alpha$. Suppose that the decomposition $A=A_0\oplus A_1$ is a
${\mathbb Z}_2$-grading such that the generators $a,b, z_1^{(1)}$ are homogeneous. Then the
graded PI-exponent $exp^{gr}(A)$ exists and
$$
exp^{gr}(A)=exp(A)=\Phi(\frac{1}{m+\alpha}),
$$
where
$$
\Phi(x)=\frac{1}{x^x(1-x)^{(1-x)}}.
$$
Moreover, if $A^\#$ is obtained from $A$ by adjoining an external unit with the induced
${\mathbb Z}_2$-grading then $exp^{gr}(A^\#)$ also exists and
$$
exp^{gr}(A^\#)=exp^{gr}(A)+1.
$$
\hfill $\Box$
\end{theorem}

\section*{Acknowledgements}
We express our sincere thanks to the referee for the numerous comments and suggestions. 
The first author was supported by the Slovenian Research Agency grants BI-RU/16-18-002 and P1-0292.
The second author was supported by the Russian Science Foundation grant 16-11-10013.


\begin{thebibliography}{99}

\bibitem{BDLaa}
Yu. Bahturin, V. Drensky,
Graded polynomial identities of matrices, {\em
Linear Algebra Appl.} {\bf 357} (2002) 15-34.

\bibitem{BBZ}
O. E. Bezushchak, A. A. Beljaev, M. V. Zaitsev. Exponents of identities of algebras with
adjoint unit, {\em Vestnik Kievskogo Nats. Univ. Ser. pfiz.-mat. nauk} 3 (2012) 7-9.

\bibitem{Dr}
V. Drensky,{\em Free algebras and PI-algebras. Graduate course in
algebra} (Springer-Verlag Singapore, 2000).

\bibitem{GMVZ}
A. Giambruno, S. P. Mishchenko, A. Valenti, M. V. Zaicev,
Polynomial codimension growth and the Specht problem, {\em  J. Algebra} {\bf 469} (2017), 421–436.

\bibitem{GMZ2}
A. Giambruno, S. P. Mishchenko, M. V. Zaicev, Algebras with intermediate growth
of the codimensions, {\em Adv. in Appl. Math.} {\bf 37}(3) (2006) 360-377.

\bibitem{GMZAdv}
A. Giambruno, S. Mishchenko, M. Zaicev, Codimensions of algebras and growth functions,
{\em Adv. Math.} {\bf 217}(3) (2008) 1027-1052.

\bibitem{G-R}
A. Giambruno, A. Regev,  Wreath products and P.I. algebras,
{\em J. Pure Appl. Algebra} {\bf 35}(2) (1985) 133-149.

\bibitem{GRZ2}
A. Giambruno, A. Regev, M. Zaicev,  On the codimension growth
of finite-dimensional Lie algebras, {\em J. Algebra} {\bf 220}(2) (1999) 466-474.

\bibitem{GRZ1}
A. Giambruno, A. Regev, M. Zaicev, Simple and semisimple
 Lie algebras and codimension growth, {\em Trans. Amer. Math. Soc.} {\bf 352}(4) (2000) 1935-1946.

\bibitem{GShZ}
A. Giambruno, I. Shestakov, M. Zaicev, Finite-dimensional
non-associative algebras and codimension growth, {\em Adv. in Appl. Math.}
{\bf 47} (2011) 125-139.

\bibitem{GZ1}
A. Giambruno, M. Zaicev,
On codimension growth of finitely generated associative algebras,
{\em Adv. Math.} {\bf 140}(2) (1998) 145-155.

\bibitem{GZ2}
A. Giambruno, M. Zaicev, Exponential codimension growth of P.I.~algebras:
an exact estimate, {\em Adv. Math.} {\bf 142}(2) (1999) 221-243.

\bibitem{GZbook}
A. Giambruno, M. Zaicev,
{\em Polynomial Identities and Asymptotic Methods}
Mathematical Surveys and Monographs, 122, Amer. Math. Soc. (Providence R.I., 2005).

\bibitem{GZ3}
A. Giambruno, M. Zaicev, Proper identities, Lie identities and exponential
codimension growth. {\em J. Algebra} {\bf 320}(5) (2008) 1933-1962.

\bibitem{GZJLMS2012}
A. Giambruno, M. Zaicev,
 On codimension growth of finite dimensional Lie superalgebras,
{\em J. London Math. Soc.} {\bf 95} (2012) 534-548.

\bibitem{GZForum}
A. Giambruno, M. Zaicev, Anomalies on codimension growth of algebras,
 {\em Forum Math.} {\bf 28}(4) (2016) 649-656.

\bibitem{L}
M. Lothaire,
{\em Algebraic Combinatorics on Words},
Encyclopedia Math. Appl., vol. 90 (Cambridge University Press, Cambridge, 2002).

\bibitem{Rats}
S. M. Ratseev, Correlation of Poisson algebras and Lie algebras in the language of identities,
Translation of {\em Mat. Zametki} {\bf 96}(4) (2014) 567-577.  {\em Math. Notes}
 {\bf 96}(3-4) (2014) 538-547.

\bibitem{RZ-CA}
D. Repov\v s, M. Zaicev,  Numerical invariants of identities of unital algebras,
{\em Comm. Algebra} {\bf 43}(9) (2015) 3823-3839.

\bibitem{Z}
M. Zaitsev,
Integrality of exponents of growth of identities of finite-dimensional Lie algebras, (Russian)
{\em Izv. Ross. Akad. Nauk Ser. Mat.} {\bf 66} (2002)(3) 23-48;
translation in {\em Izv. Math.} {\bf 66}(3) (2002) 463-487.

\bibitem{Z-AL}
M. V.  Zaitsev,  Identities of finite-dimensional unitary algebras,(Russian)
{\em  Algebra  Logika} {\bf 50}(5) (2011) 563-594, 693, 695; translation in
{\em  Algebra Logic} {\bf 50}(5) (2011) 381-404.

\bibitem{ERA}
M. Zaicev, On existence of PI-exponents of codimension growth, {\em Electron. Res.
  Announc. Math. Sci.} {\bf 21} (2014) 113-119.

\bibitem{VMGU}
M. V. Zaicev, Graded identities in finite-dimensional algebras of codimensions of identities
in associative algebras, Translation of {\em  Vestnik Moskov. Univ. Ser. I Mat. Mekh.} (5) 2015 54-57.
{\em Moscow Univ. Math. Bull.} {\bf 70}(5) (2015) 234-236.

\bibitem{RZ-MatSb}
M. V. Zaicev, D. Repov\v s, Exponential growth of codimensions of identities of algebras
with unity, (Russian) {\em Mat. Sb.} {\bf 206}(10) (2015) 103-126; translation in
{\em  Sb. Math.} {\bf 206}(10) (2015) 1440-1462.

\end{thebibliography}
\end{document}